\RequirePackage{fix-cm}
\documentclass[smallextended]{svjour3}       
\smartqed  
\usepackage{graphicx}
\usepackage{lipsum}
\usepackage{amsfonts}
\usepackage{graphicx}
\usepackage{amssymb}
\usepackage{mathtools} 
\usepackage{url}
\usepackage{caption}
\usepackage{subfig}

\newcommand{\R}{\ensuremath{\mathbb{R}}}

\begin{document}

\title{Data Driven Conditional Optimal Transport}

\author{Esteban G. Tabak  \and Giulio Trigila \and Wenjun Zhao}

\institute{Esteban G. Tabak \at
              Courant Institute of Mathematical Sciences, 251 Mercer Street, New York, NY 10012, USA \\
              \email{tabak@cims.nyu.edu}           
              \and
              Giulio Trigila \at
              Baruch College, CUNY, 55 Lexington avenue, New York, NY 10035, USA \\
              \email{giulio.trigila@baruch.cuny.edu}
              \and
             Wenjun Zhao\at
             Courant Institute of Mathematical Sciences, 251 Mercer Street, New York, NY 10012, USA \\
             \email{wenjun@cims.nyu.edu}
}

\date{Received: date / Accepted: date}

\maketitle

\begin{abstract}
A data driven procedure is developed to compute the optimal map between two conditional probabilities $\rho(x|z_{1},...,z_{L})$ and $\mu(y|z_{1},...,z_{L})$ depending on a set of covariates $z_{i}$. The procedure is tested on synthetic data from the ACIC Data Analysis Challenge 2017 and it is applied to non uniform lightness transfer between images. Exactly solvable examples and simulations are performed to highlight  the differences with ordinary optimal transport.

\end{abstract}

\section{Introduction}
Optimal transport seeks the mass preserving map $T$ between two probability distributions that minimizes the expected value of a given cost function, the \emph{transportation cost} between a point and its image under $T$ \cite{Vil}. The minimal cost defines a metric in the space of probability distributions, the \emph{Wasserstein distance}. Beyond providing a metric, the optimal map $T$ itself has broad applicability, which this article extends through the development of conditional optimal transport.

Consider as a specific example the evaluation of the effects of a long-term medical treatment (alternatively of a habit, such as smoking or dieting). Optimal transport can be used to quantify the changes in probability distribution of quantities that  characterize the health state of a person (blood pressure, blood sugar level, heart beat rate) in the two scenarios: with and without treatment. Data typically consist of independent measurements of these quantities in treated and untreated populations. Yet the distribution of these quantities depends on many covariates beyond the presence or absence of treatment, such as age, weight, sex, habits. Hence one seeks the effect of the treatment as a function of these covariates.

Motivated by this and similar applications, this article develops a data driven procedure to compute the optimal map $T(x, z)$ between two conditional probability densities $\rho(x|z_{1},...,z_{L})$ and $\mu(y|z_{1},...,z_{L})$, with covariates $z_{i}$. In the example above, $y = T(x, z)$ estimates the value $y$ that the quantity of interest would have under treatment if, without treatment, its value were $x$, under specific values of the covariates $z_l$. The procedure is data driven, as it uses only samples $\left\{x^i, z_1^i, \ldots, z_L^i \right\}$ and $\left\{y^i, z_1^i, \ldots, z_L^i \right\}$  from $\rho$ and $\mu$. Notice that we do not seek a pairwise matching between $\left\{x^i, z_1^i, \ldots, z_L^i \right\}$ and $\left\{y^i, z_1^i, \ldots, z_L^i \right\}$: typically these two data sets do not even have the same cardinality. Instead, we  work under the hypothesis that these samples are drawn from smooth conditional densities $\rho(x|z)=\rho(x,z)/\gamma^x(z)$, $\mu(y|z)=\mu(y,z)/\gamma^y(z)$ and covariate distributions $\gamma^x(z)$ and $\gamma^y(z)$, and hence we seek a map $y = T(x, z)$ that is a smooth function of its arguments.

The need for conditional optimal transport is particularly apparent when the distributions for the covariates $z$ for the source and target distributions are unbalanced, i.e. when $\gamma^x$ and $\gamma^y$ are different. Consider as a particularly telling example a situation when the treatment has no effect, i.e. $\rho(x|z) = \mu(x|z)$, so we should have $y = x$, yet $\gamma^x \ne \gamma^y$:
$$ \rho(x|z) = \mu(y|z) = N(z, 1), \quad \gamma^x(z) = N(-1, 1), \quad \gamma^y(z) = N(1, 1),$$
where $N(a,b)$ denotes the $1$d normal distribution with mean $a$ and variance $b$.
Then 
$$ \rho(x) = \int \rho(x|z) \gamma^x(z) \ dz = N(-1, 2), \quad \mu(y) = \int \mu(y|z) \gamma^y(z) \ dz = N(1, 2).$$
It follows that, if one would not look at the covariate $z$, one would infer incorrectly that $y = x+2$, i.e. that the treatment does have a significant effect. We will see in section \ref{sec:treatment} an instance of this phenomenon appearing in the more complex setting of a biomedical application, where conditional transport provides critical aid.

\begin{figure}[!htb]
  \begin{center}
      \begin{tabular}{ccc}                                              
      \resizebox{45mm}{!}{\includegraphics{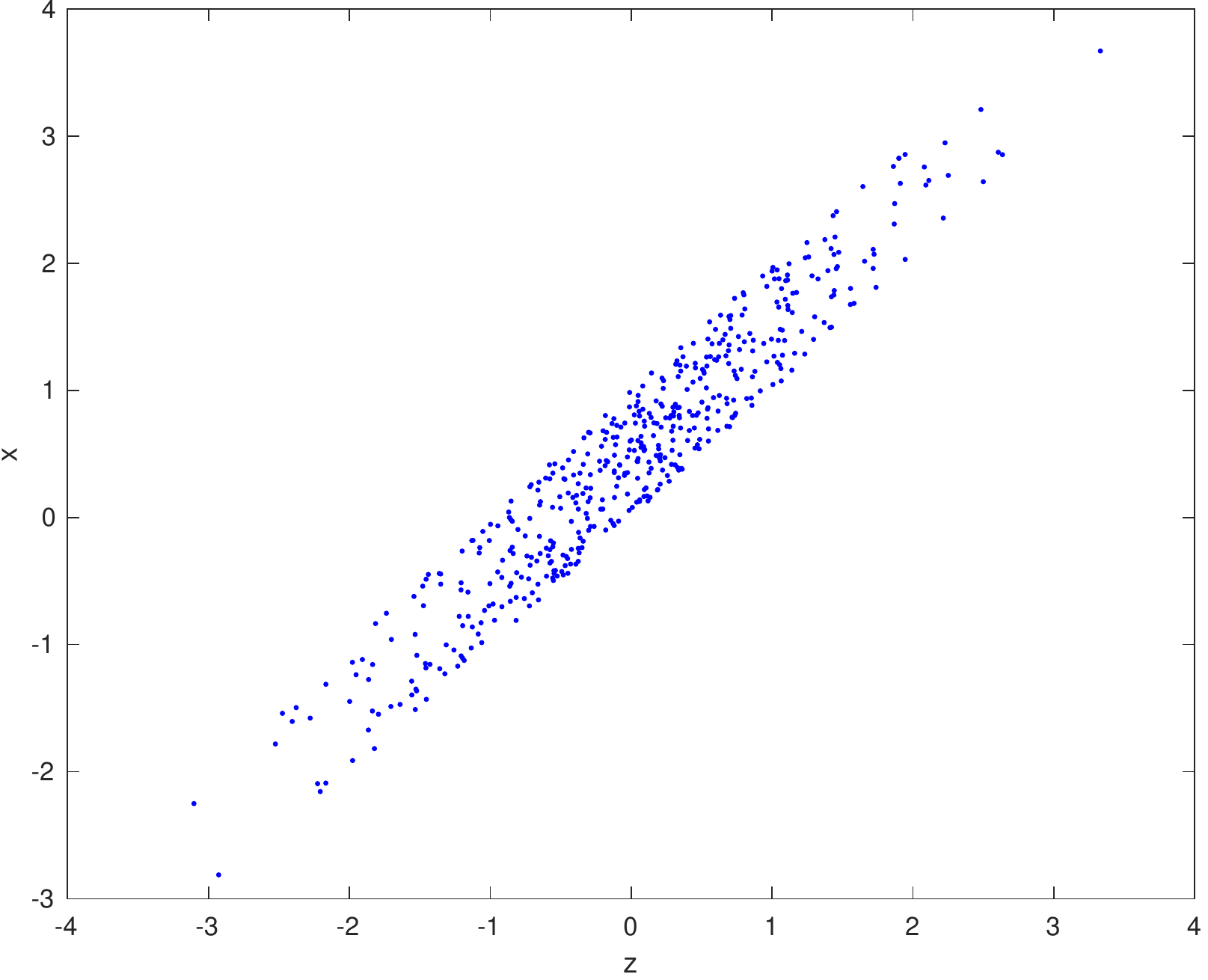}} &        
            \resizebox{45mm}{!}{\includegraphics{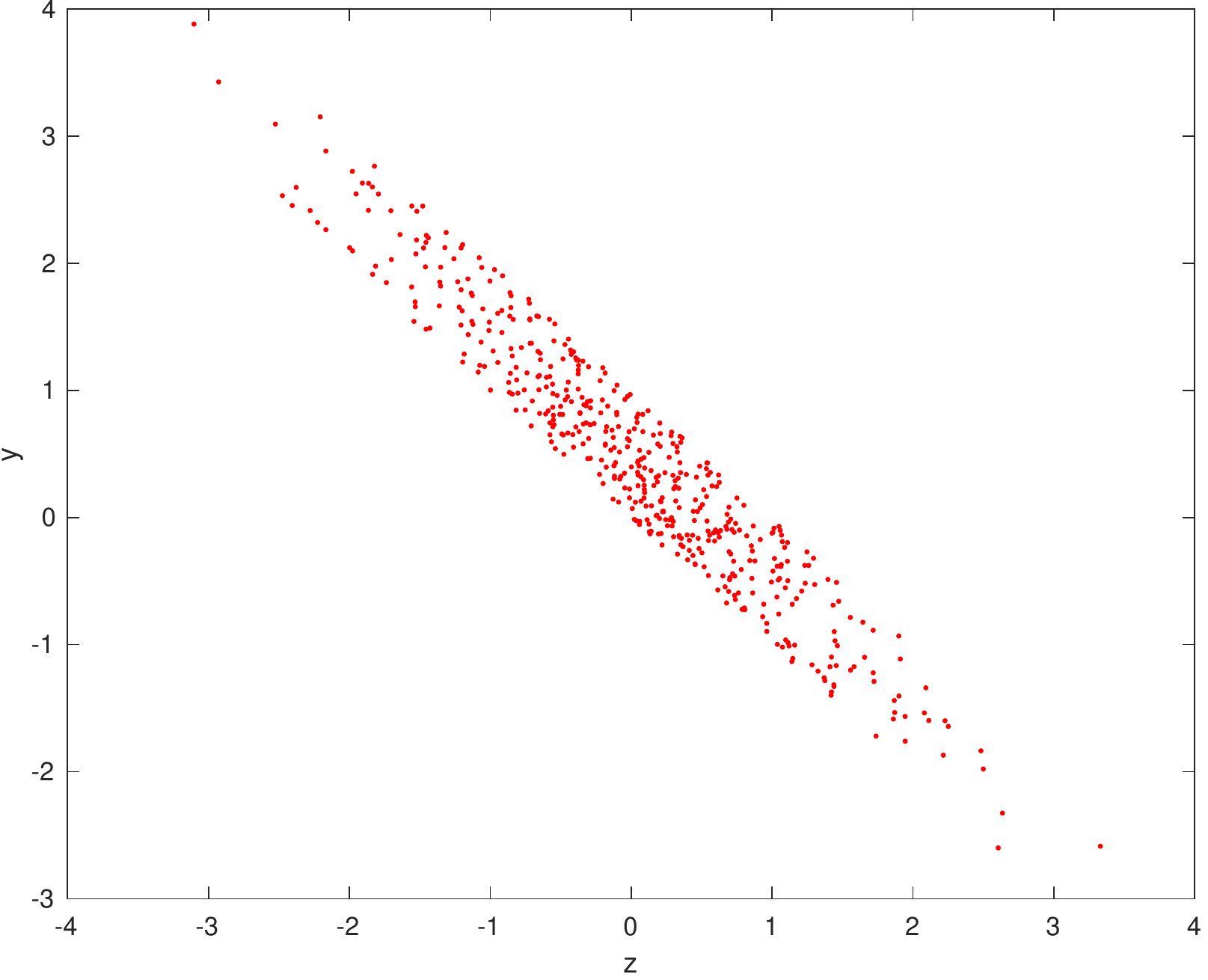}}\\                           
	\resizebox{45mm}{!}{\includegraphics{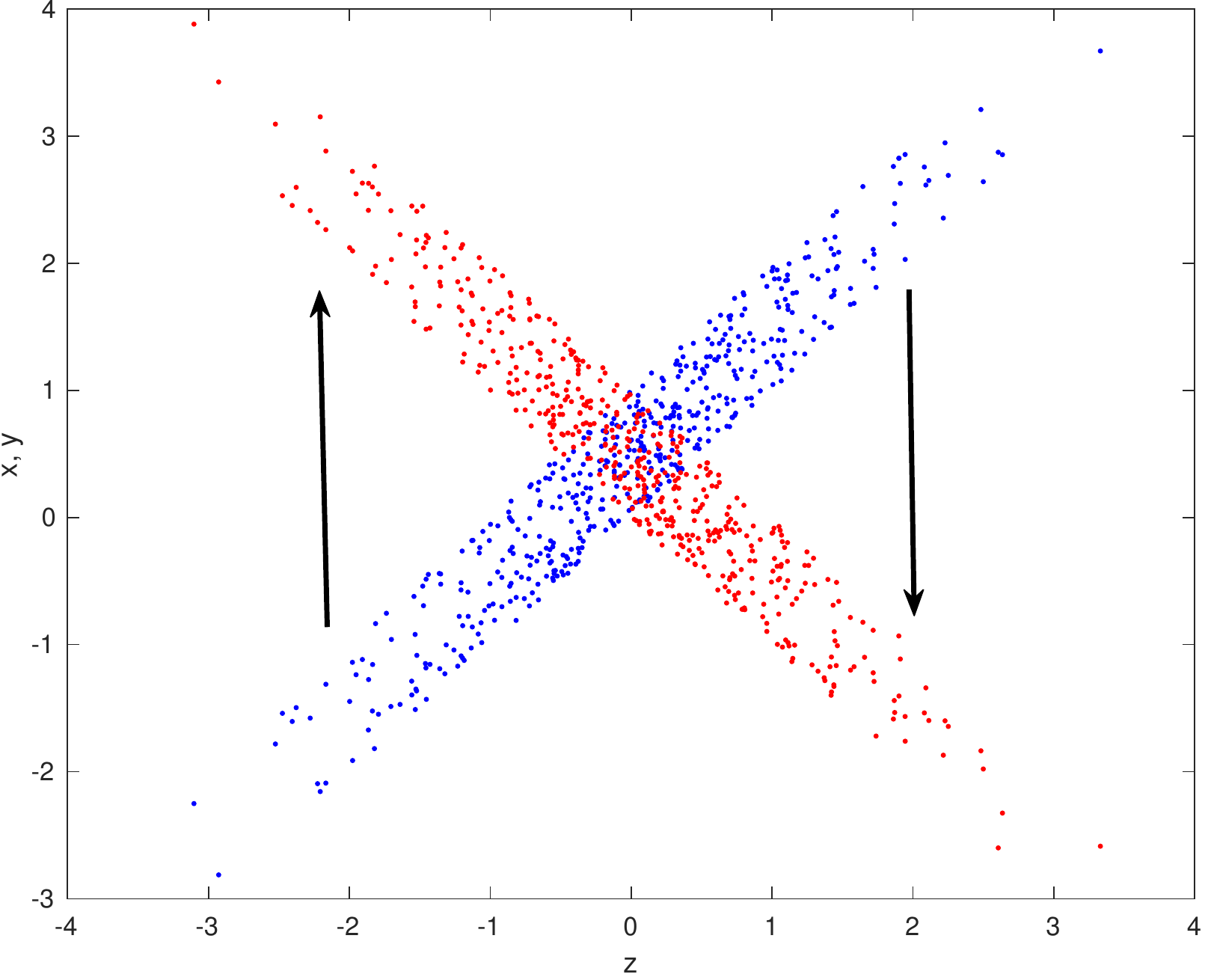}}&
      \resizebox{45mm}{!}{\includegraphics{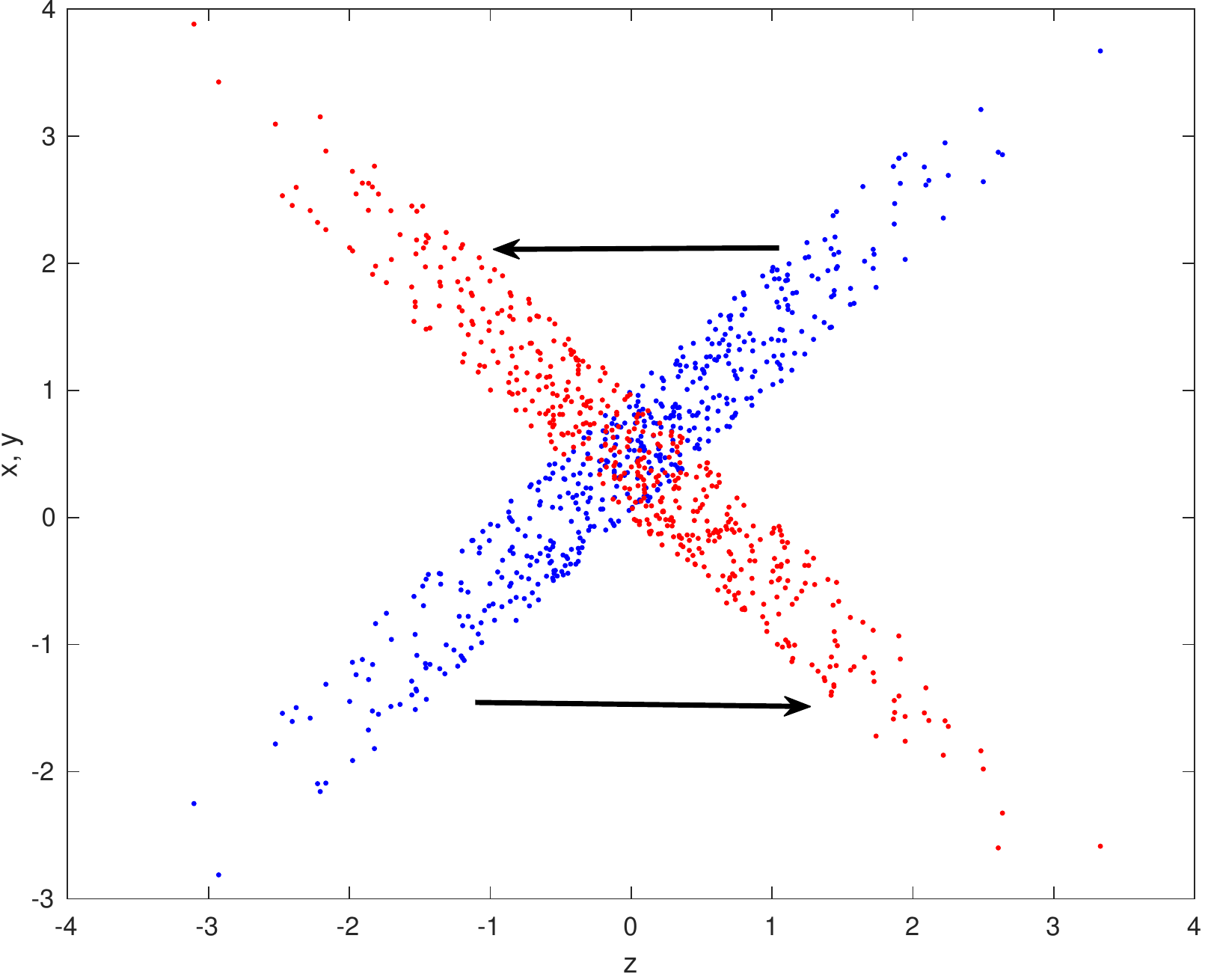}}
     \end{tabular}
    \caption{Upper row: source (left) and target (right) distributions. Lower left: optimal transport of $x$ conditioned on $z$, the arrows indicate that the lower left branch and the upper right branch of the source distribution are mapped respectively to the upper left branch and the lower right branch of the target distribution. Lower right: optimal transport of $z$ conditioned on $x$, in this case is the upper right branch of the source distribution to be mapped to the upper left branch of the target distribution.}
    \label{fig:syntexample}
  \end{center}
\end{figure}

Conditional transport provides a very flexible toolbox for data analysis, as the choice of which variables are conditioned to which others is left at the discretion of the analyst. In anticipation of the application of this principle to color transfer problems in section \ref{sec:lightness}, we illustrate it here with a simple example. Consider a covariate $z \sim N(0, 1)$ and two dependent variables $x \sim N(z, 1)$ and $y \sim N(-z, 1)$ (see Figure \ref{fig:syntexample} for a sketch relative to this problem). Since the marginals $\rho(x)$ and $\mu(y)$ are identical, performing optimal transport between them yields the identity map $y = x$, while conditioning to $z$ yields $y = x - 2 z$,  effectively rotating the joint distribution $\rho(x, z)$ clockwise, and performing two dimensional transport between $\rho(x, z)$ and $\mu(y, z)$ yields an irrotational map \cite{Vil}. Finally, if in a thought experiment we would identify $x$ and $y$ and switch the roles of dependent and independent variables, conditioning the transport in $z$-space to $x$, we would obtain $z_2 = z_1 - 2x$, effectively rotating the joint distribution $\rho(x, z)$ counter-clockwise.

\section{Conditional optimal transport}
\label{sec:COT}

Conditional optimal transport between two conditional distributions $\rho(x|z)$ and $\mu(y|z)$ can be defined simply as the map $T(x, z)$ that performs optimal transport between them for each value of $z$:
\begin{equation}
\label{eq:OT0}
\forall z \ \begin{dcases}
\min_{T(:, z)}\int c(T(x, z),x)\rho(x|z) dx\\
T\#\rho(:|z)=\mu(:|z),
\end{dcases}
\end{equation}
where $c(x,y)$ represents the cost of moving a unit of mass from $x$ to $y$  and the symbol $\#$ indicates the push forward of probability measures, i.e. if $x$ has distribution $\rho(x|z)$ then $y = T(x, z)$ has distribution $\mu(y|z) = T\#\rho(:|z)$. Since $T(:, z)$ decouples under different values of $z$, we can multiply the cost by the distribution $\gamma^x(z)$ of the covariates $z$ in the source and integrate over $z$, yielding
\begin{equation}
\label{eq:OT1}
\begin{dcases}
\min_{T(:, z)}\int c(T(x, z),x)\ \rho(x, z)\ dx dz\\
\forall z\ T\#\rho(:|z)=\mu(:|z),
\end{dcases}
\end{equation}
where $\rho(x|z)$ denotes the conditional and $\rho(x, z) = \rho(x|z) \gamma^x(z)$ the joint distribution of $x$ and $z$.

We need to reformulate this problem in a way that is implementable in terms of samples $\left\{x^i, z^i\right\}$ and $\left\{y^j, z^j\right\}$. As is stands in (\ref{eq:OT1}), two immediate problems emerge: there are not enough samples for each value of $z$, typically none or  one for continuous covariates, to characterize the corresponding conditional distributions, and it is not clear how to enforce or verify the push forward condition. The first problem is at the very heart of the need for conditional optimal transport: even though the objective functions for each value of $z$ decouple, one thinks of a commonality across $z$ that makes samples from each conditional distribution be informative on the others. In the case of continuous covariates $z$, this can be posed as a smoothness condition (in $z$) on $\rho(x|z)$.

In order to address the second problem, we interpret the push forward condition in terms of relative entropy:
$$ T\#\rho(:|z)=\mu(:|z) \iff D_{KL}(T\#\rho(:|z), \mu(:|z)) = 0, $$
where
$$
D_{KL}(\rho_1(x|z) || \rho_2(x|z) ) = \int \gamma_1(z) \int \log \left( \frac{\rho_1(x|z)}{\rho_2(x|z)} \right)\rho_1(x|z) \ dx dz
$$
is the conditional Kullback-Leibler divergence between $\rho_1$ and $\rho_2$ (\cite{cover2012elements}). Since this is non-negative, we can rewrite the problem in (\ref{eq:OT1}) as
$$
\min_{T(: ,z)} \max_{\lambda \ge 0} \left[ \int c(T(x, z),x)\rho(x, z) \ dx dz + \lambda D_{KL}(T\#\rho(:|z), \mu(:|z)) \right] . 
$$
Instead of maximizing over $\lambda$, it will be convenient to fix a value of $\lambda$ large enough that the push forward condition can be considered satisfied for all practical purposes (it is straightforward to prove that, as $\lambda \rightarrow \infty$ the solution with fixed $\lambda$ converges to the true minimax solution. In our implementation below, $\lambda$ grows at each step of the algorithm.) Then the problem above becomes
$$
\min_{T} \left[ \int c(T(x, z),x)\rho(x, z) \ dx dz + \lambda\ D_{KL}(T\#\rho(:|z), \mu(:|z)) \right], \quad \lambda \gg 1.
$$
For any  $\rho_1(x, z) = \gamma_1(z)  \rho_1(x|z)$ and $\rho_2(x, z) = \gamma_2(z)  \rho_2(x|z)$, we have the ``chain rule'' for the relative entropy (\cite{cover2012elements}),
$$
D_{KL}(\rho_1(x|z) || \rho_2(x|z) ) = D_{KL}(\rho_1(x,z) || \rho_2(x,z) ) - D_{KL}(\gamma_1(z) || \gamma_2(z) ).
$$
Since the map $T$ acts only on $x$, it has no effect the last term, so we can write
$$
\min_{T} \left[ \int c(T(x, z),x)\rho(x, z) \ dx dz + \lambda D_{KL}(T\#\rho(x, z), \mu(x, z)) \right], \quad \lambda \gg 1.
$$
This formulation improves over the one in (\ref{eq:OT0}) by consolidating an infinite set of problems, one for every value of $z$, into a single one. Yet it is not clear yet how to enforce the push forward condition in terms of samples, as the definition of the relative entropy involves logarithms of $\rho$ and $\mu$. To address this, we invoke a variational formulation of the relative entropy between two distributions \cite{donsker1975asymptotic}:
\begin{equation}
D_{KL}(\rho, \mu) =
\max_g \left[ \int g(x) \rho(x) dx - \log\left(\int e^{g(x)} \mu(x)\ dx\right) \right],
\end{equation}
which involves $\rho$ and $\mu$ only in the calculation of the expected values of $g$ and $e^g$, with a natural sample-based interpretation as empirical means. Then our problem becomes
\begin{equation}
 \label{theor_COT}
\min_{T} \max_g \int c(T(x, z),x) d\rho(x, z) 
 + \lambda \left[\int g(T(x, z), z) d\rho(x, z) - \log\left(\int e^{g(y, z)} d\mu(y, z)\right) \right] 
\end{equation}
or, in terms of samples,
\begin{equation}
\min_{T} \max_g \Bigg[\frac{1}{n} \sum_i \Big(c(T(x_i, z_i), x_i)\rho(x_i, z_i)
+ \lambda g(T(x_i, z_i), z_i)\Big)
 - \lambda \log\left(\frac{1}{m} \sum_j e^{g(y_j, z_j)}  \right) \Bigg].
\label{sample_COT}
\end{equation}
This adversarial formulation has two players with strategies $T$ and $g$, one minimizing the cost and the other enforcing the push forward condition, providing an adaptive ``lens'' that identifies those places where the push-forward condition does not hold: for any $T$, the optimal $g$ in (\ref{theor_COT}) is given by
$$ g = \log\left(\frac{\rho(T(x,z)|z)}{\mu(y|z)}\right) + \log\left(\frac{\gamma^x(z)}{\gamma^y(z)}\right) ,$$
where the first term is furthest from zero in those places where $\rho(T(x,z)|z)$ and $\mu(y|z)$ differ the most.

\section{Parametrization of the flows}
\label{sec:Flow}
In order to complete the problem formulation in (\ref{sample_COT}), we need to specify the family of functions over which the map $T(x, z)$ and the test-function $g(y, z)$ are optimized. These families should satisfy some general properties:
\begin{enumerate}

\item  be rich enough that $g$ can capture all significant differences between $\rho(x|z)$ and $\mu(y|z)$ and $T$ can resolve them,

\item not be so rich as to overfit the sample points $\left\{x^i, z^i\right\}$,  $\left\{y^j, z^j\right\}$. For instance, a $g$ with arbitrarily small bandwidth would force the sets $\left\{T(x^i, z^i), z^i\right\}$, $\left\{y^j, z^j\right\}$ to agree point-wise, an extreme case of overfitting that is not only undesirable but also unattainable when their cardinality differs. Moreover, the dependence of the functions on $z$ should be such that, with a finite number of samples, it should still capture the assumed smoothness of $\rho(x|z)$: functions  that are too localized in $z$ space effectively decouple the transport problems for every value of $z$, for which there are not enough available sample points,

\item be well-balanced: if one of the two players has a much richer toolbox than the other, the game would be unfair, leading not only to a waste of computational resources but also possibly to instability and inaccuracy, and

\item be apt to robust and effective optimization.  
\end{enumerate}

These conditions leave space for many proposals, such as defining $T$ and $g$ through neural networks. Instead, the examples in this article are solved with the two approaches detailed below. Both share the feature that $T$ is built on map composition: at each step $n$ of the mini-maximization algorithm, an elementary map $E^n$ is applied not to the original sample points $\left\{x^i\right\}$, but to their current images:
%
$$ T^n(x^i, z^i) = E^n\left(T^{n-1}(x^i, z^i), z^i\right). $$
This way, simple elementary maps $E$ depending on only a handful of parameters can give rise through map composition to rich global maps $T$. The two proposals differ in that one builds nonlinear richness through evolving Gaussian mixtures, while the other builds complex $z$-dependence through an extra compositional step. In this article, the first method is applied to a lightness transfer problem, and the second to the effect of a medical treatment, as the latter is linear in $x$ but has complex, nonlinear dependence on many covariates $z$. 
\subsection{Evolving Gaussian mixtures}
We adopt as elementary map the gradient of a convex potential function: $E(x, z) = \nabla_x \Phi(x,z)$, with $\Phi$ built from a quadratic form in $x$ with coefficients that depend on $z$, plus a combination of Gaussians in $(x, z)$ space, and similarly for the test function $g$. By having the centers and amplitudes of these Gaussians evolve, we can approximate quite general functions $\Phi$ and $g$.

Notice that the gradient of a radial basis function kernel with bandwidth $d$,
$$
G_d(x,x') = \exp \left( - \frac{||x - x'||^2}{2d^2}\right),
$$
is bounded by $ \pm \frac{1}{d\exp(1/2)}$, and its second order derivatives by $ \frac{2}{d^2 \exp(3/2)}<\frac{1}{2d^2}$. It follows that $ \frac{1}{2d^2}\frac{||\mathbf x||_2^2}{2} \pm G_d (\textbf z,\textbf m_{z_i})G_d(\textbf x,\textbf m_i)$ is convex, so we propose
\begin{multline}
\Phi(\mathbf x,\mathbf z) = (\mathbf c^T_0  + \mathbf z^T \mathbf c_1 )  x + 
\frac 12 \mathbf x^T \mathbf C_2(\mathbf z)\mathbf x 
+ \sum_{i=1}^K a_i^2   \left(\frac{||\mathbf x||_2^2}{4d^2}  - G_{d} (\mathbf z,\mathbf  m_{z_i}) G_d (\mathbf x,\mathbf m_i) \right) + \\
 \sum_i b^2_i  \left(  \frac{||\mathbf x||_2^2}{4d^2}  + G_{d} (\mathbf z, \mathbf  m_{z_i}) G_{d} (\mathbf x, \mathbf m_i)   \right), \
\mathbf C_2(\mathbf z) = \mathbf C_{2,0} ^T\mathbf C_{2,0} + \mathbf z^T  \mathbf  C_{2,1}^T \mathbf  C_{2,1} \mathbf z,
\nonumber
\end{multline}
with $\mathbf C_{2,0}, \mathbf C_{2,1}$ lower triangular.
In order to start the map at every step at the identity, the initialization must satisfy
$$
\mathbf C_{2,0}(i,i)^2 + \sum_i^K \frac{1}{4d^2}( a_i^2 + b_i^2) = 1, \quad a^2_{i} =b_i^2,
$$
so we propose
$$
a^2_{i} =b_i^2 =  \frac{4d^2 \delta}{ 2 K  },\quad \mathbf C_{2,0}(i,i) = \sqrt{1-\delta}, \quad \delta = \frac{1}{2},
$$
with all other parameters starting from zero.
The bandwidth $d$ is chosen via $d = quantile(pdist([y; z]),1/K)$,
where $pdist$ is the pairwise distance function. With this choice there are approximately $1/K$ points in the effective support of each Gaussian.

For the test function, we propose 
$$
g(\mathbf x,\mathbf z) = \sum_{i=1}^K \alpha_i  G_{d} (\mathbf z ,\mathbf m_{z_i}) G_d( \mathbf x, \mathbf  m_i) + ( \mathbf \beta_0 ^ T  +  \mathbf z ^T \mathbf \beta_1 ) \mathbf x + \mathbf x^T ( \mathbf \beta_2+ \sum_i \mathbf \beta_{3,i} z_i )  \mathbf x,
$$
with each iteration starting at the parameter values from the previous step.
The Gaussian centers are treated differently in the test function $g$, where they are extra parameters to ascend, and in the potential $\Phi$, where they are fixed at their values from $g$ in the prior step. The underlying notion is that $g$ locates those areas where the distributions do not agree, and then $T$ corrects them.
\subsection{Extended map composition}
This second methodology considers maps given by rigid translations and test functions that capture the conditional mean $\bar{x}(z)$:
$$ T(x, z) = x + U(z), \quad g(y, z) = V(z) y + W(z), \quad x \in \R,$$
with general, nonlinear dependence on $z$. To build these, we define a composition function 
$$
F(\mathbf a, z, v, u)  =  (a^1_0 + \sum_{i=1}^L a^1_i z_i + a^1_{L+1}u )  + (a^2_0 + \sum_{i=1}^L a^2_i z_i + a^2_{L+1}u ) v,
$$
in terms of which the test function at each step is given by
$$
g^{n+1}(y, z) = v^{n+1} y  + w^{n+1}  ,\quad v^{n+1} = F(\mathbf \beta, z, v^n, u^n), \quad w^{n+1} = F(\mathbf \eta, z, w^n, 0),
$$
and the map by
$$
T ^{n+1} (T^n, z) =  T^n + u^{n+1}, \quad u^{n+1} =F(\mathbf \alpha, z, u^n, v^n).
$$
These maps are initialized at $u^0=v^0=w^0=0$. Before each each step, $\alpha$ is set to $0$ (as $T$ is reinitialized every step to the identity), and so are $\beta$ and $\eta$, except for $\beta_0^2=\eta_0^2 = 1$, which makes $g$ evolve from its value at the previous step. 

\section{Examples}\label{sec:examples}
We illustrate the procedure with two applications: determination of the effect of a medical treatment and lightness transfer. In order to solve the problem (\ref{sample_COT}) we use the general procedure for mini-maximization described in \cite{Minimax}.  
\subsection{Effect of a Treatment}\label{sec:treatment}
\begin{figure}[!htb]
  \begin{center}
      \begin{tabular}{ccc}  
     \includegraphics[width=43mm,height=45mm]{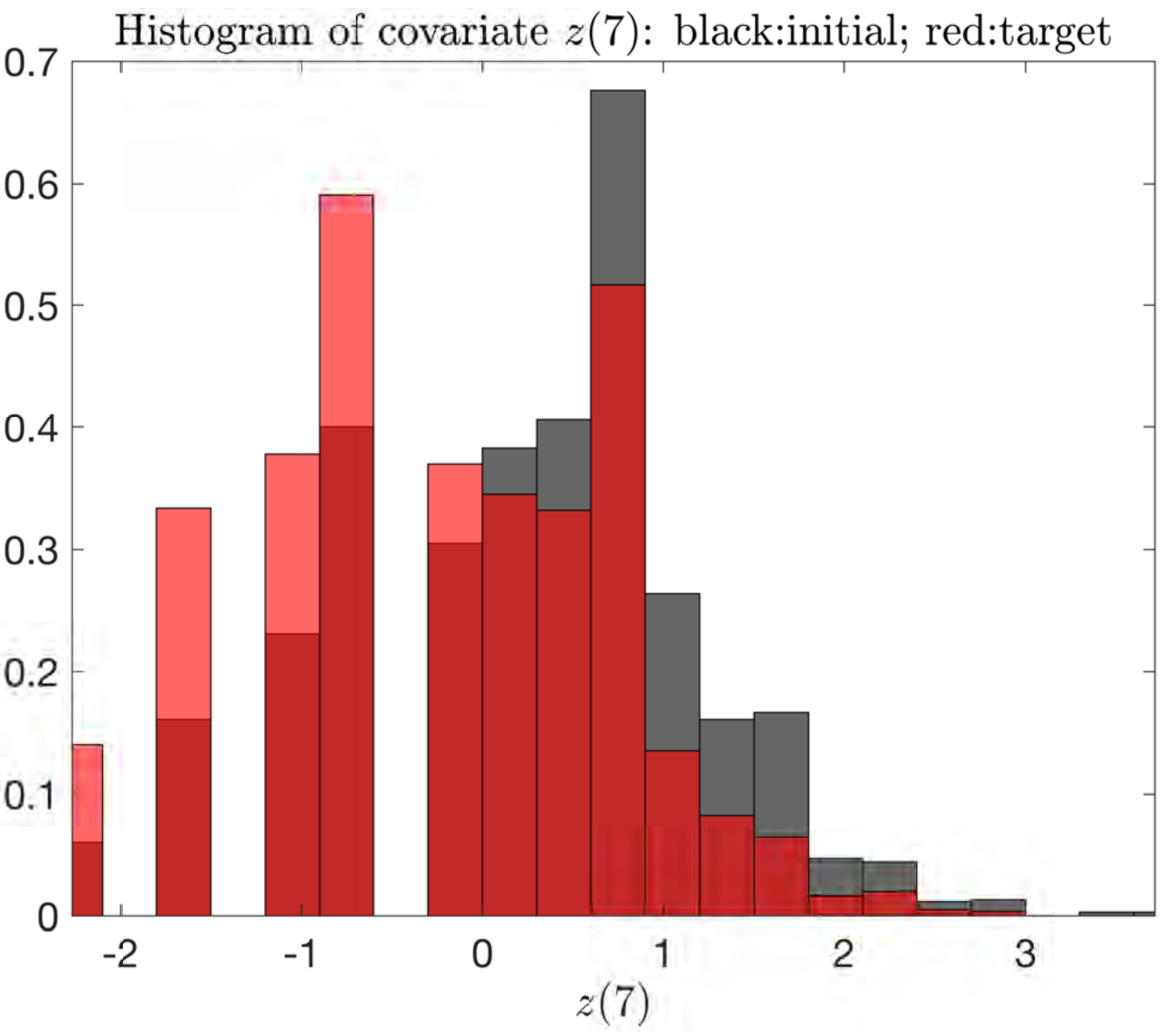} &
      \resizebox{40mm}{!}{\includegraphics{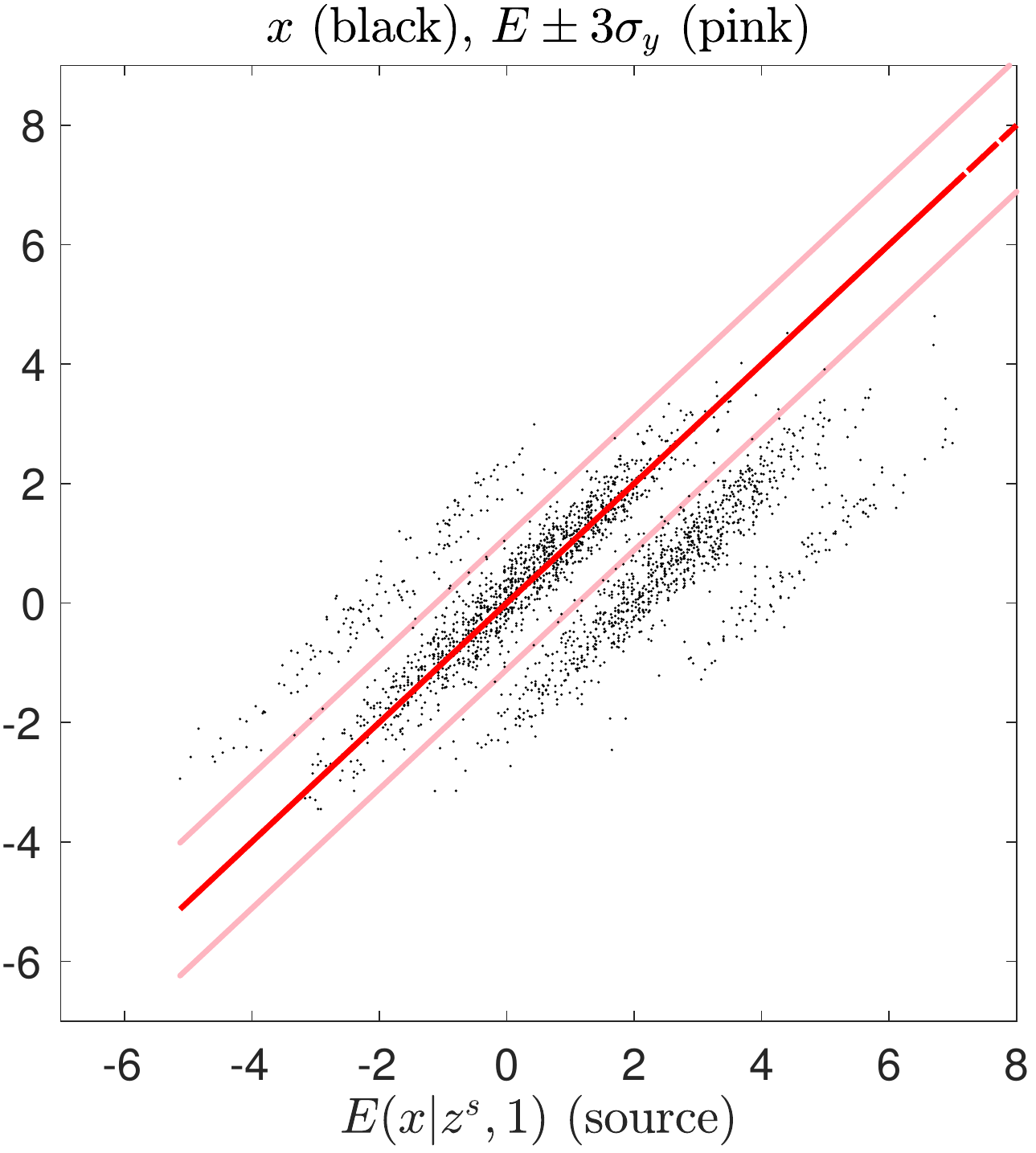}} &        
      \resizebox{40mm}{!}{\includegraphics{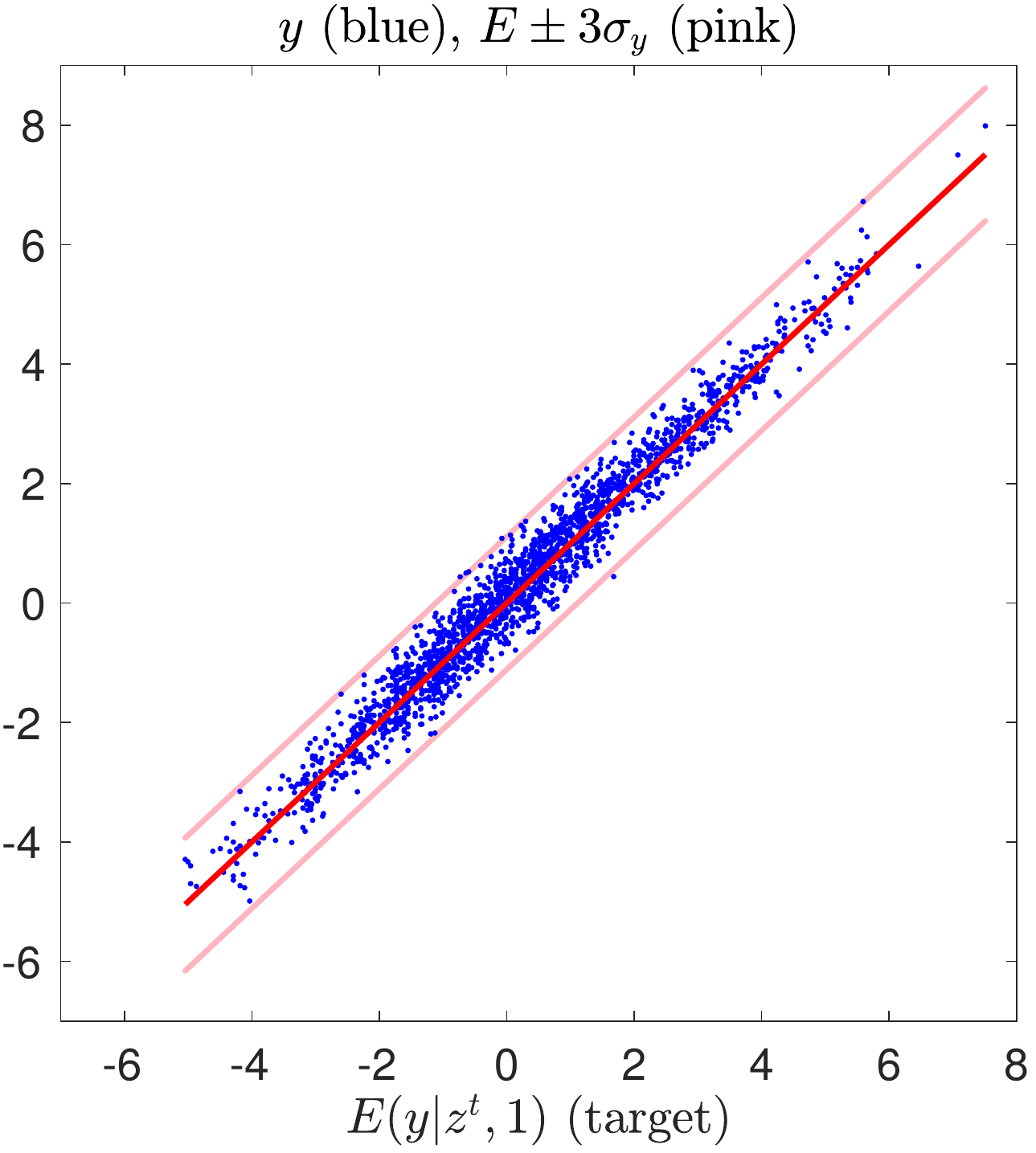}}\\                           
     \end{tabular}
    \caption{Left panel: Unbalance in the distribution of $z_{7}$ between the source and the target data set. Center: Response variable $x$ for  patients before the treatment VS theoretical expected value of the same patients undergone the treatment. Right: Response variable $y$ of treated patients VS theoretical expected value of the same patients.}
    \label{fig:Treat1}
  \end{center}
\end{figure}

\begin{figure}[!htb]
  \begin{center}
      \begin{tabular}{ccc}  
      \resizebox{40mm}{!}{\includegraphics{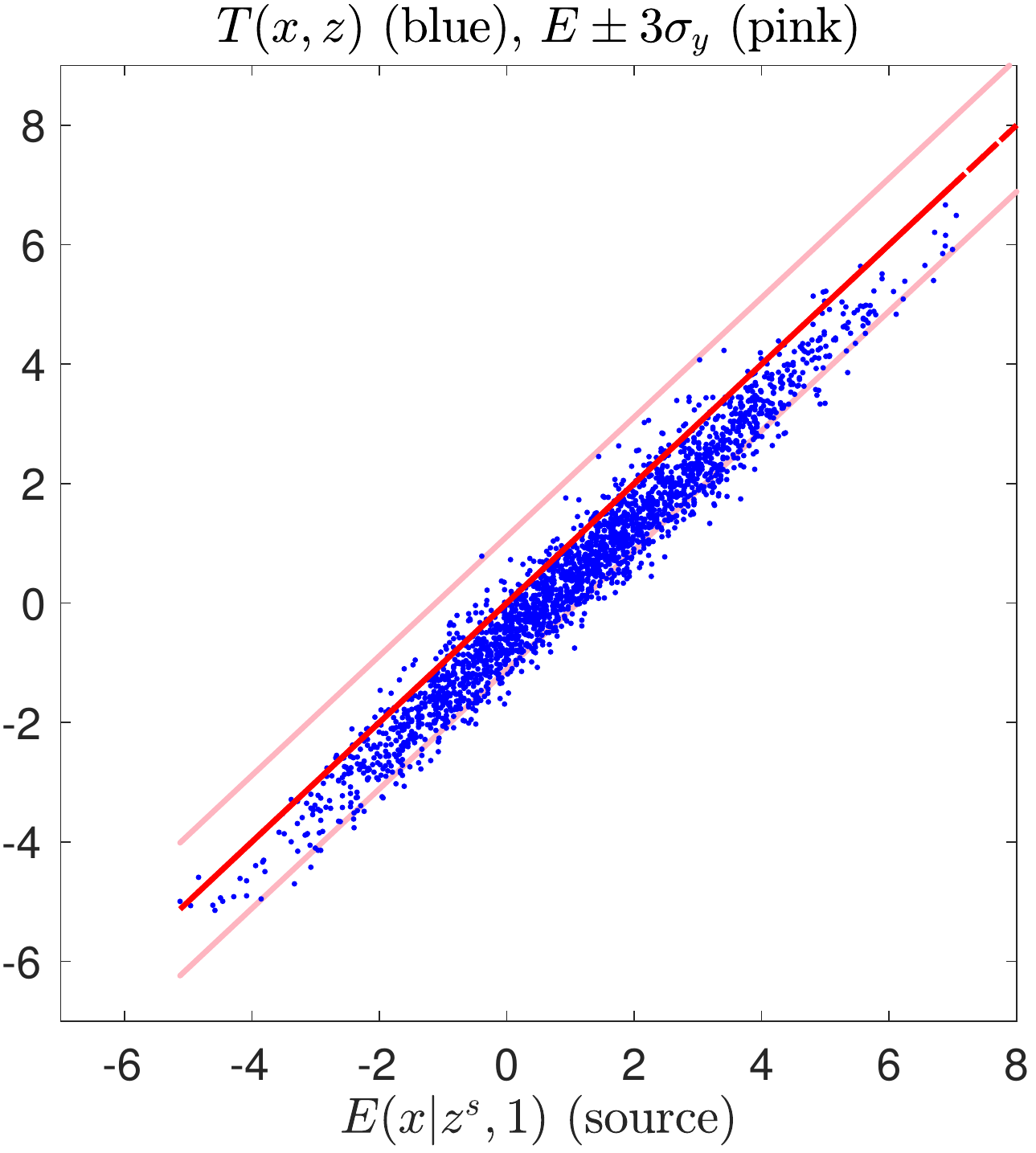}}&                                            
      \resizebox{40mm}{!}{\includegraphics{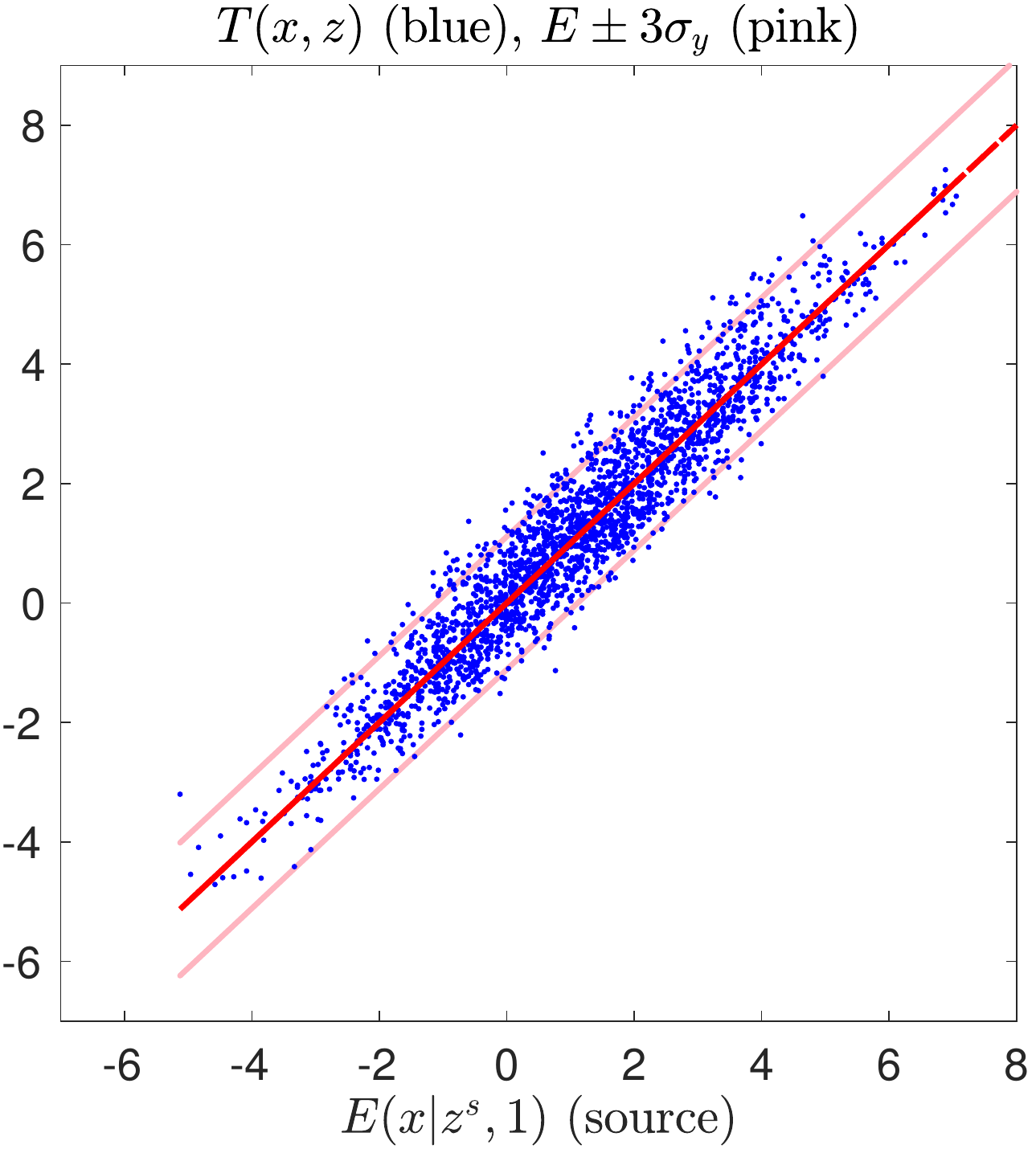}}&
      \includegraphics[width=43mm,height=45mm]{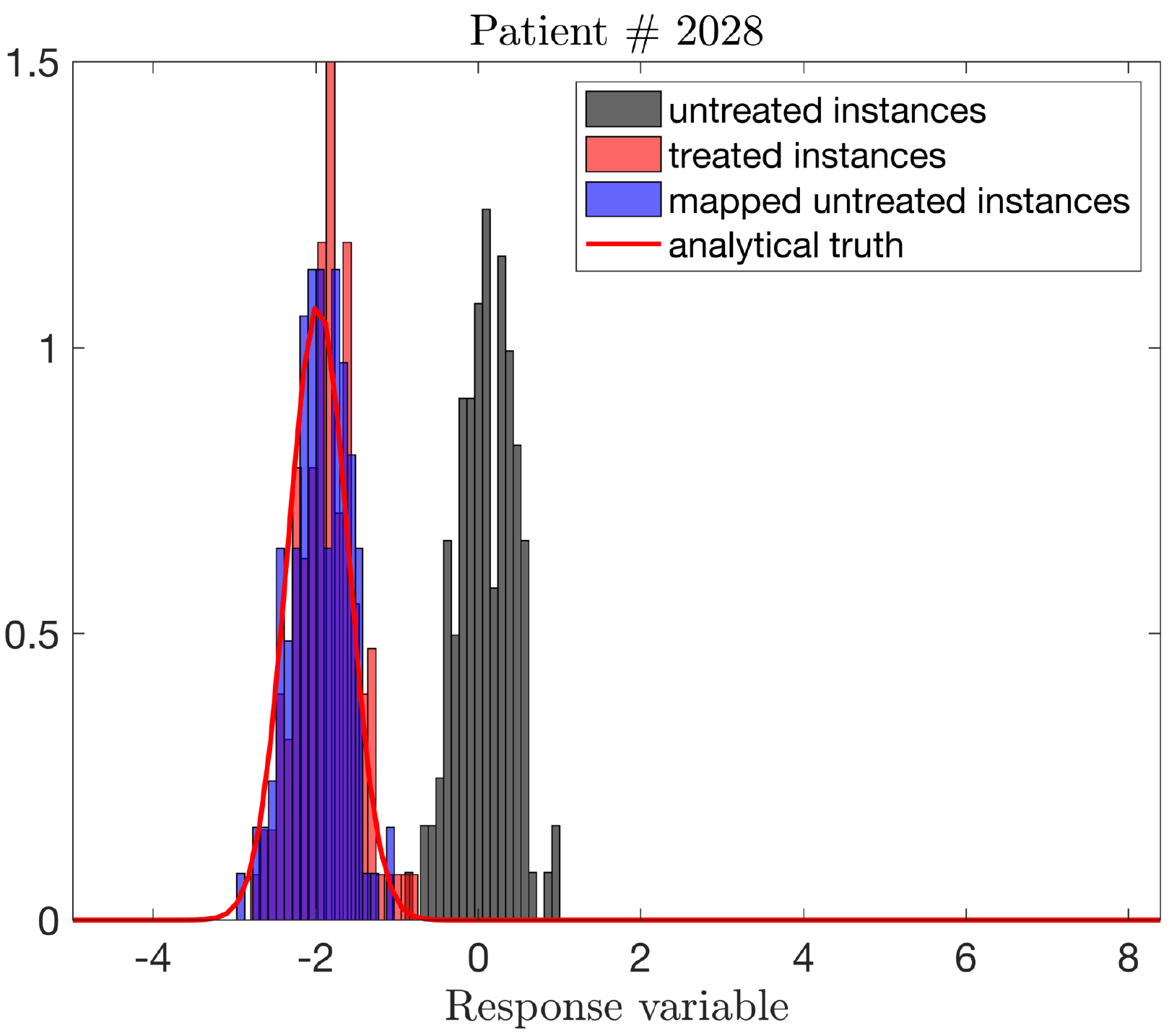}           
     \end{tabular}
    \caption{Left: numerical value of map $T(x_i, z_i)$ obtained using only the discrete covariates, which are the ones that the  true $T$ depends on. The result is biased due to the unbalance between $\gamma^x$ and $\gamma^y$ for $\gamma(z_7)$). Middle: numerical value of map $T(x_i, z_i)$ obtained using all the covariates. Right: comparison between the application of the map $T(x, z)$ to all untreated instances in the full 250 batches to the histogram of the response $y$ for all treated instances of the patient.}
    \label{fig:Treat2}
  \end{center}
\end{figure}

We apply conditional optimal transport to determine the response to a treatment of a variable $x \in \R$ in terms of covariates $z$. We use data from the ACIC data analysis challenge 2017 \cite{hahn2019atlantic} (\url{https://arxiv.org/pdf/1905.09515.pdf}), considering the first of their 32 generating models, with 8 covariates: 6 binary and 2 continuous. We divide the data set into two groups: the untreated ($x$) and treated ($y$) patients, with samples drawn from distributions $\rho(x, z) = \gamma^{x}(z) \rho(x|z)$ and $\mu(y, z) = \gamma^y(z) \mu(y|z)$, having the property that
$$ \mu(y|z) = \rho(y-\tau(z)|z), \quad \gamma^x(z) \ne \gamma^y(z). $$
It will be important for the analysis below that $\tau$ (the ``effect'' of the treatment) depends only on the binary covariates, but the marginals $\gamma(z)$ depend also on the continuous ones. We compute the optimal map $T(x, z)$ using only  the first of the 250 batches of data provided, each referring to the same 4302 patients, i.e. the same values of $z_i$ under different realizations of the noise. The middle panel of Figure \ref{fig:Treat1} displays the untreated values $x_i$ as a function of the expected value that they would have under treatment given the values ${z^s}_i$ of their covariates:
$$ E(x|z^s, 1) = \int \left(x + \tau(z^s)\right) \rho(x|z^s) \ dx, $$
while the right panel displays similarly the treated values $y_i$. The left panel of Figure \ref{fig:Treat2} displays the map $T(x_i, z_i)$ obtained using only the discrete covariates, which are the ones that the true $T$ depends on.  However, because of the unbalance between $\gamma^x$ and $\gamma^y$ (see the left panel of Figure \ref{fig:Treat1} for $\gamma(z_7)$), the results are biased, much as in the synthetic example in the introduction. The middle panel shows that, when all covariates are considered, this biased is resolved. For a specific patient, the right panel compares the application of the map $T(x, z)$ to all untreated instances in the full 250 batches to the histogram of the response $y$ for all treated instances of the patient. 
\subsection{Lightness  transfer}\label{sec:lightness}
\begin{figure}[!htb]
\begin{tabular}{cc}
	\begin{minipage}{0.5\textwidth} 
	\includegraphics[width = 0.7\textwidth]{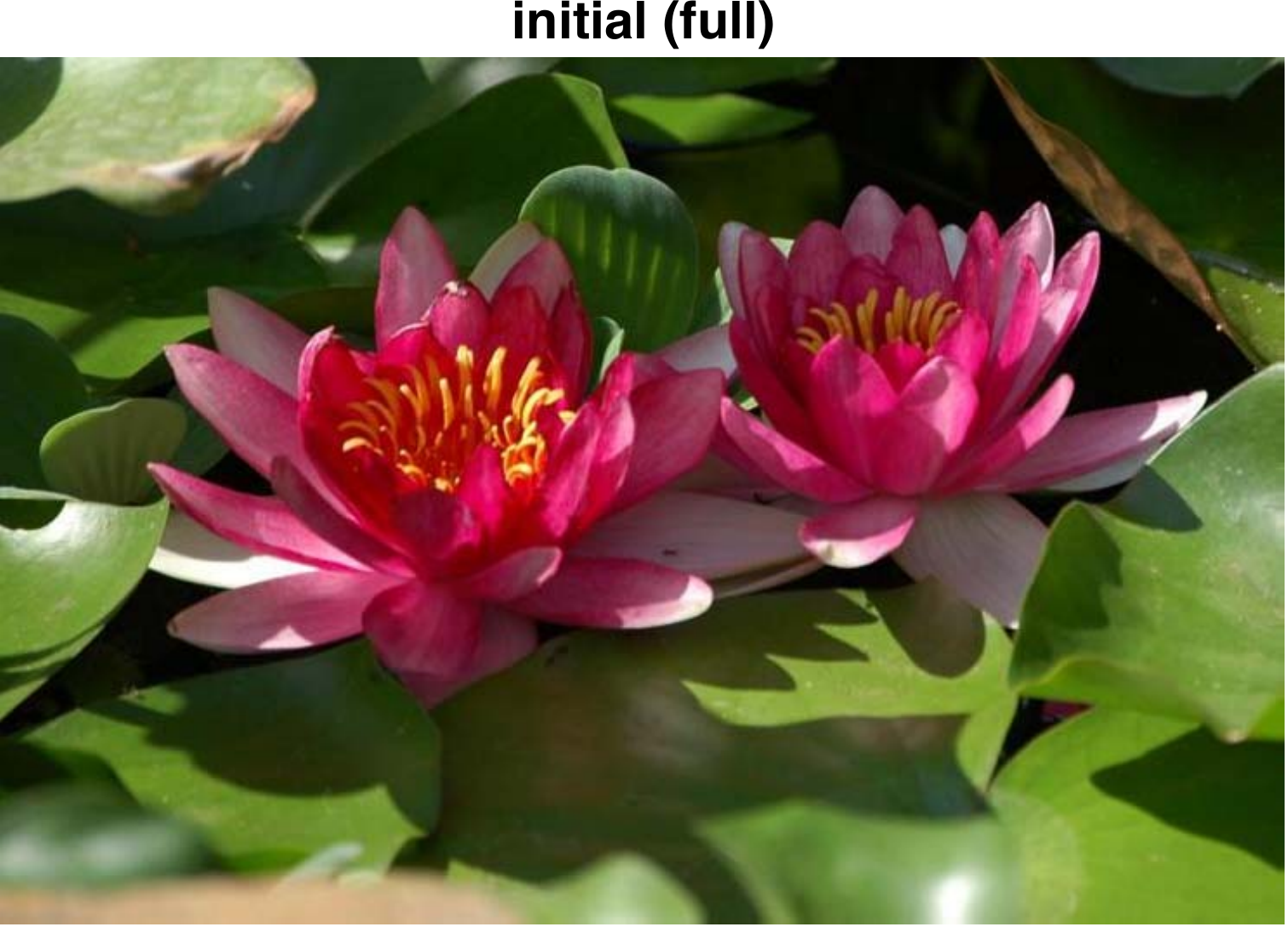} \\ 
	\includegraphics[width = 0.7\textwidth]{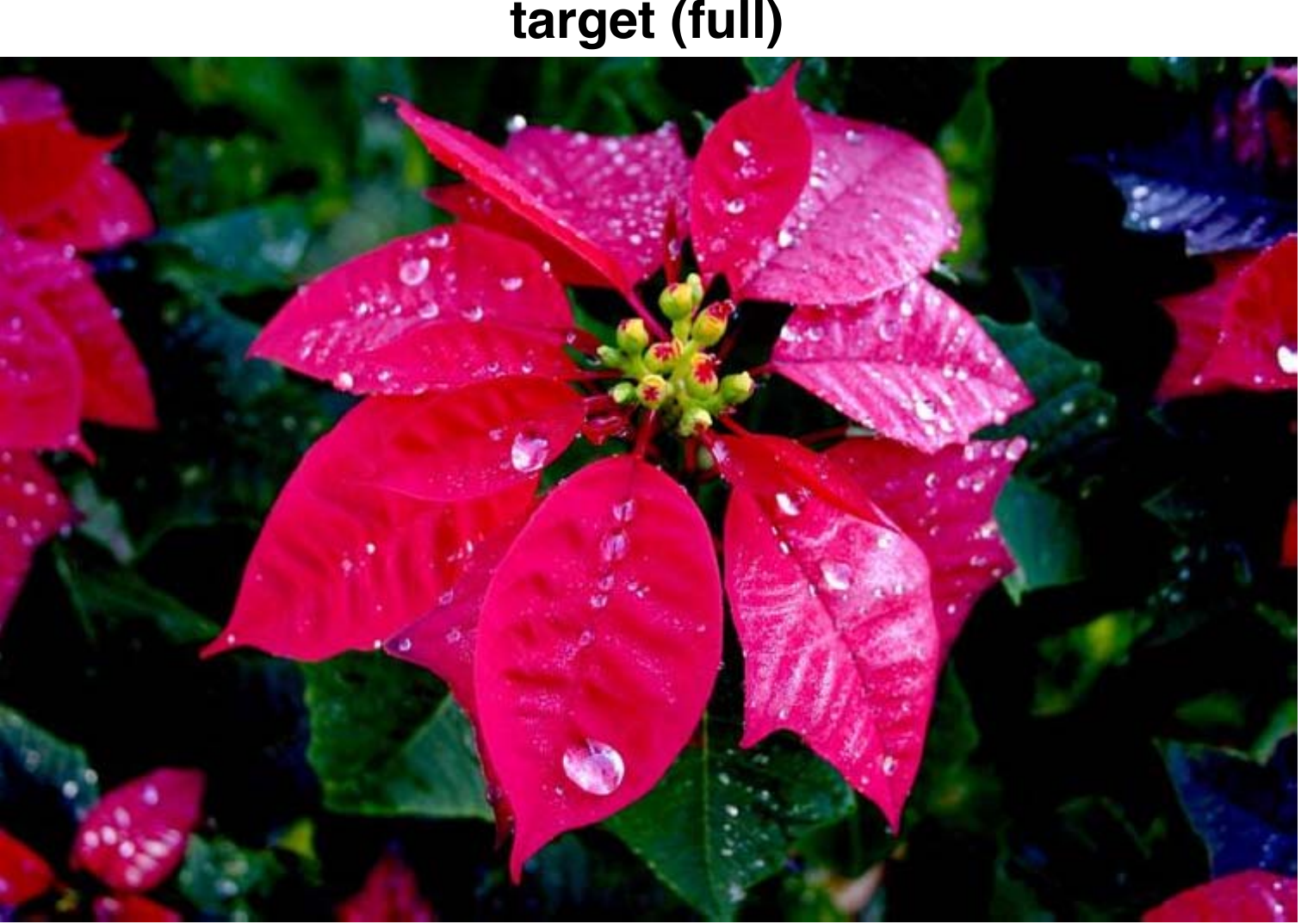}
 	\end{minipage}
 	\begin{minipage}{0.5\textwidth} 
	\includegraphics[width = 0.5\textwidth]{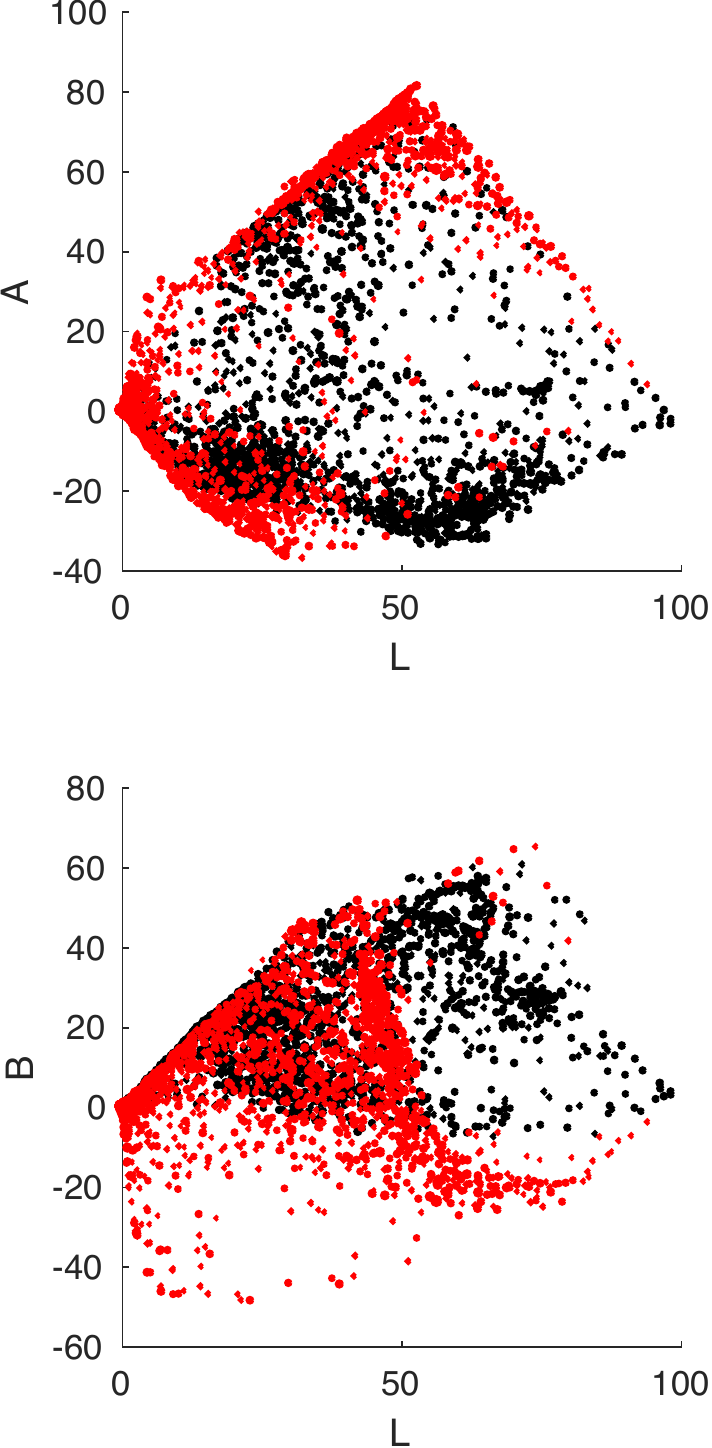} 
	\end{minipage}& 
\end{tabular}
 \caption{\label{fig:StrTarg} Left column: initial (top) and target (bottom) image. Right column: L*a*b coordinates for the initial (in red) and the target (in black) image}
 \end{figure}

Next we apply conditional optimal transport to lightness transfer. Consider the first column of Figure \ref{fig:StrTarg}, corresponding to two flowers photographed under different light conditions. We seek to transform the first photograph so as to present it under the light conditions of the second. This goes beyond merely changing lightness uniformly, since for instance at sunset certain colors are perceived as having become darker than others.

An image can be represented in the three dimensional CIELAB (L*a*b) space whose coordinates are the lightness $L$,  the red/green contrast $A$ and  blue/yellow contrast $B$. The right column of Figure \ref{fig:StrTarg} shows the images of the flowers in this L*a*b space, where each point corresponds to a superpixel, defined through a clustering procedure to introduce information about the geometry of the image \cite{rabin2014adaptive}. We follow \cite{tai2005local} to define a similarity metric by means of Gaussian kernel, map the obtained superpixels with our procedure, and use a TMR filter after the map to recover sharp details \cite{rabin2011removing}.

\begin{figure}
\subfloat[1D OT]{\includegraphics[width = 0.35\textwidth]{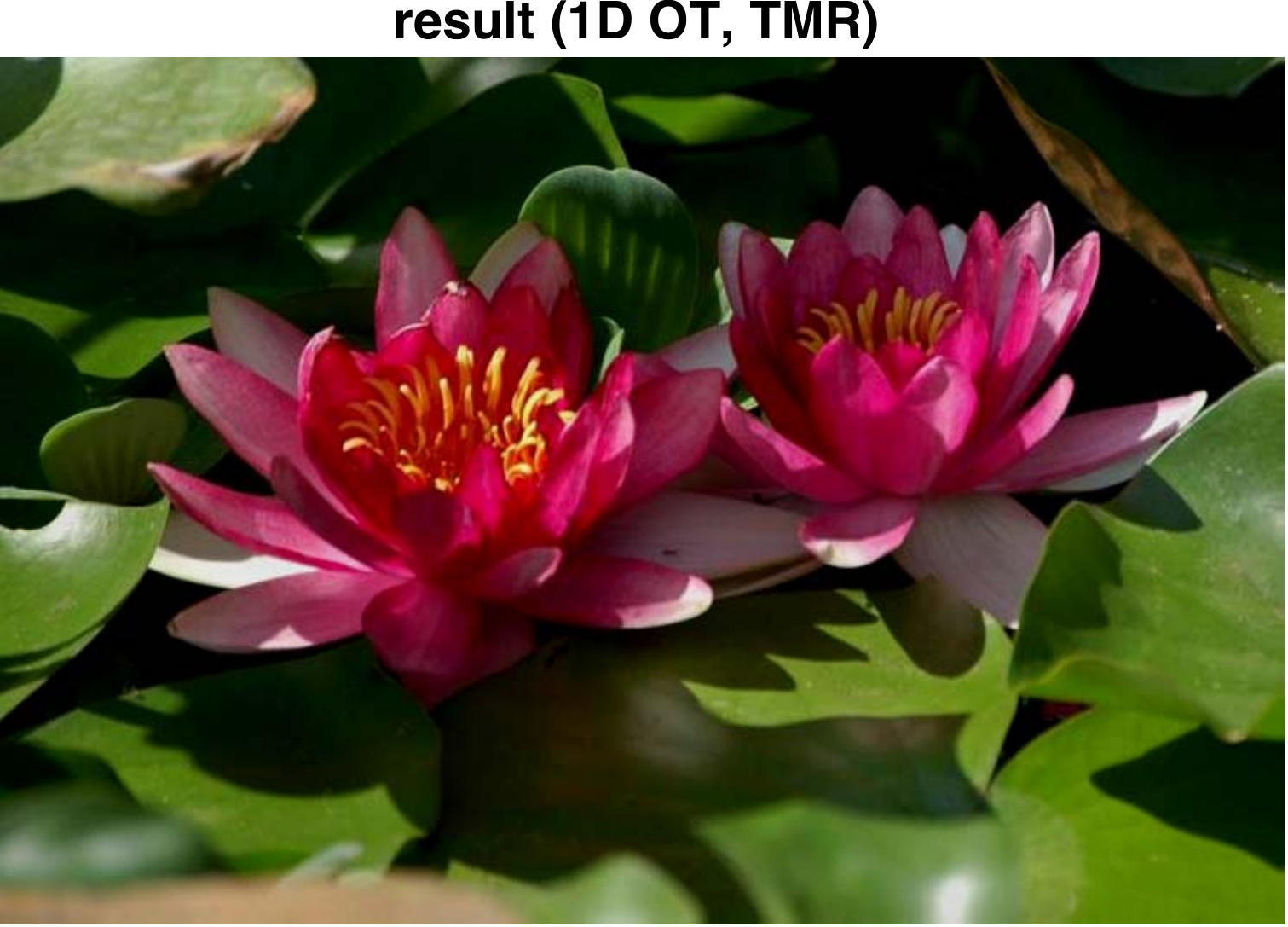}}
\subfloat[COT]{\includegraphics[width = 0.35\textwidth]{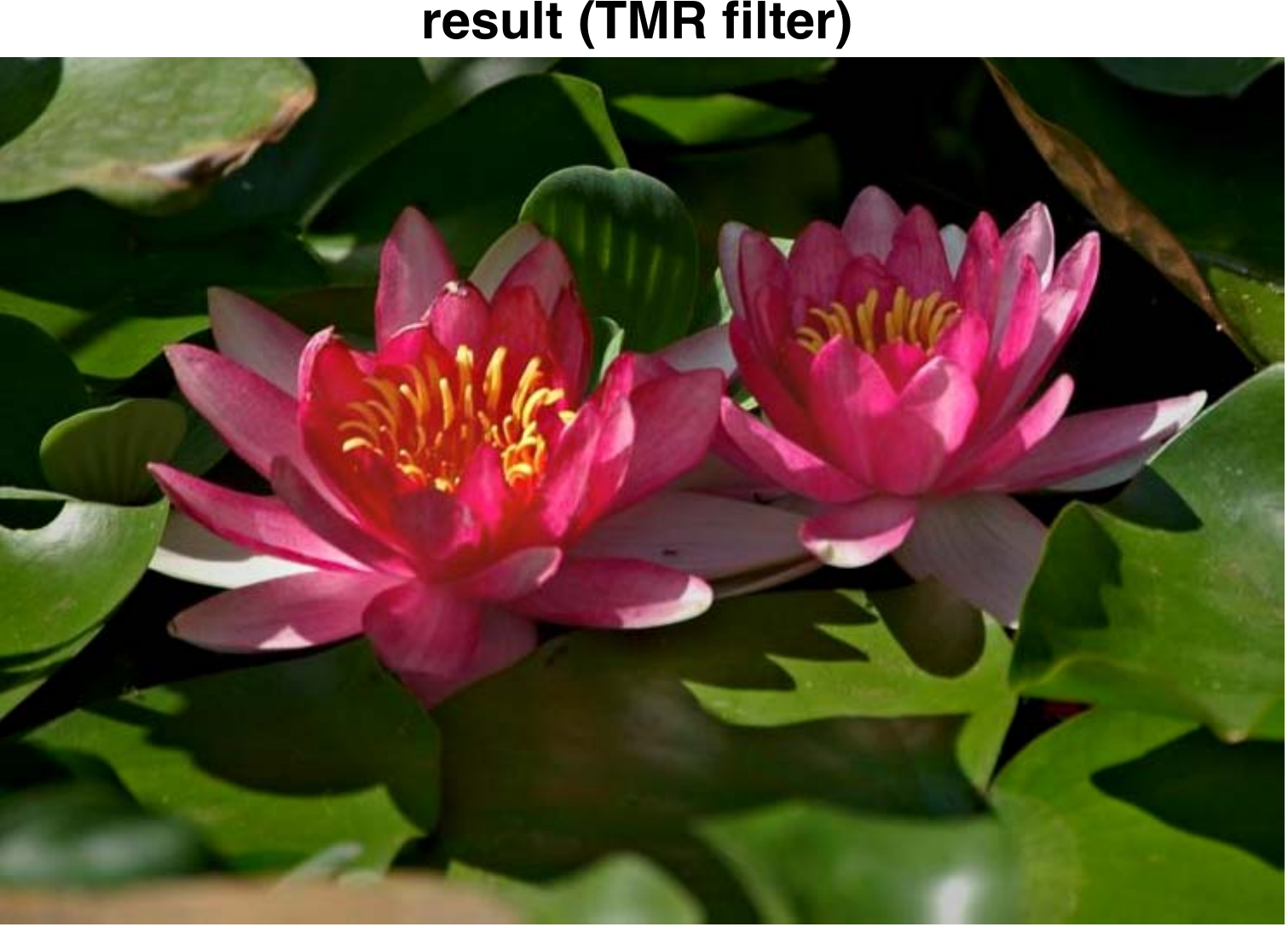}}
\subfloat[3DOT ]{\includegraphics[width = 0.35\textwidth]{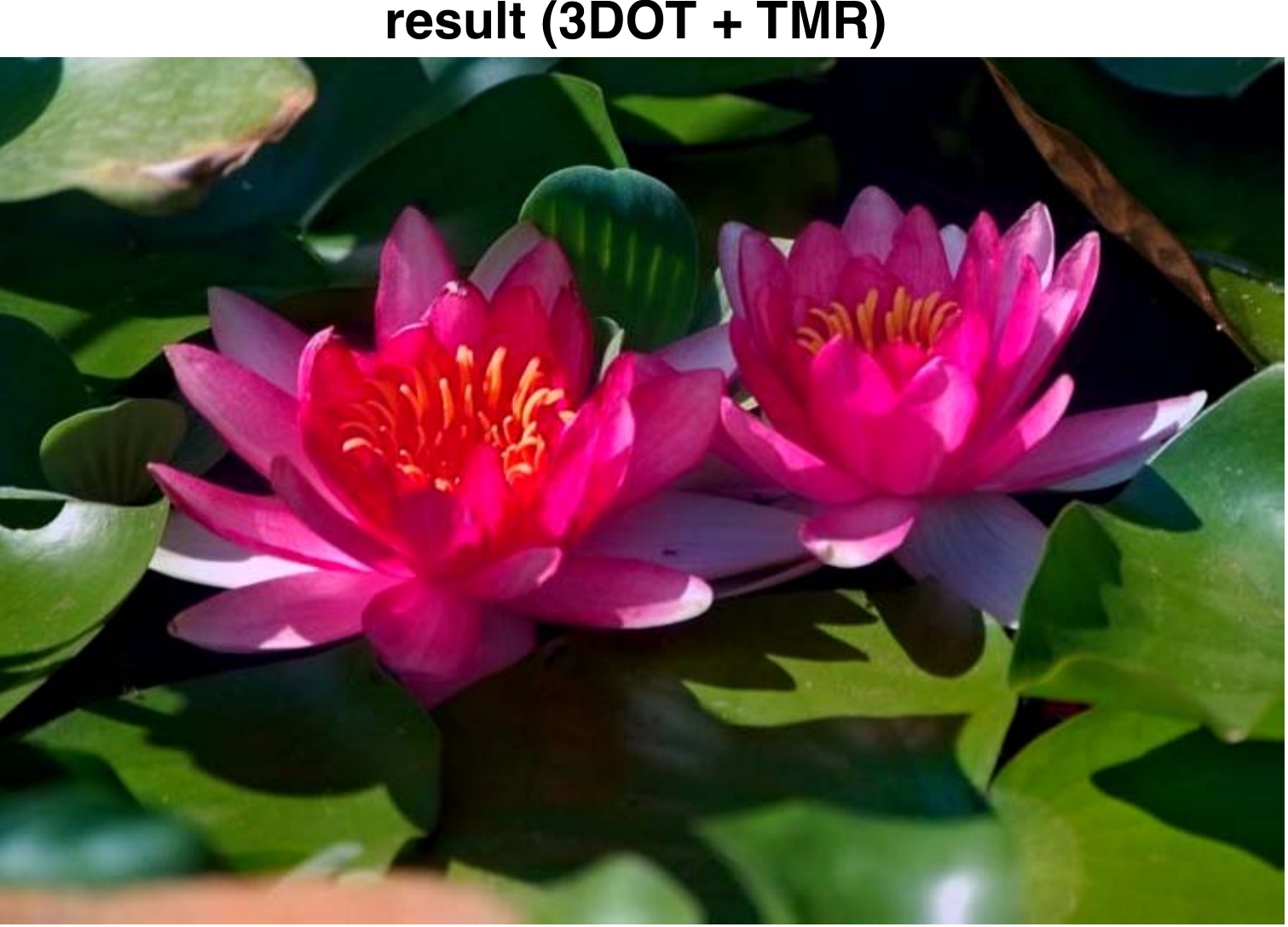}}

\subfloat[1D OT result in LAB]{\includegraphics[width = 0.35\textwidth]{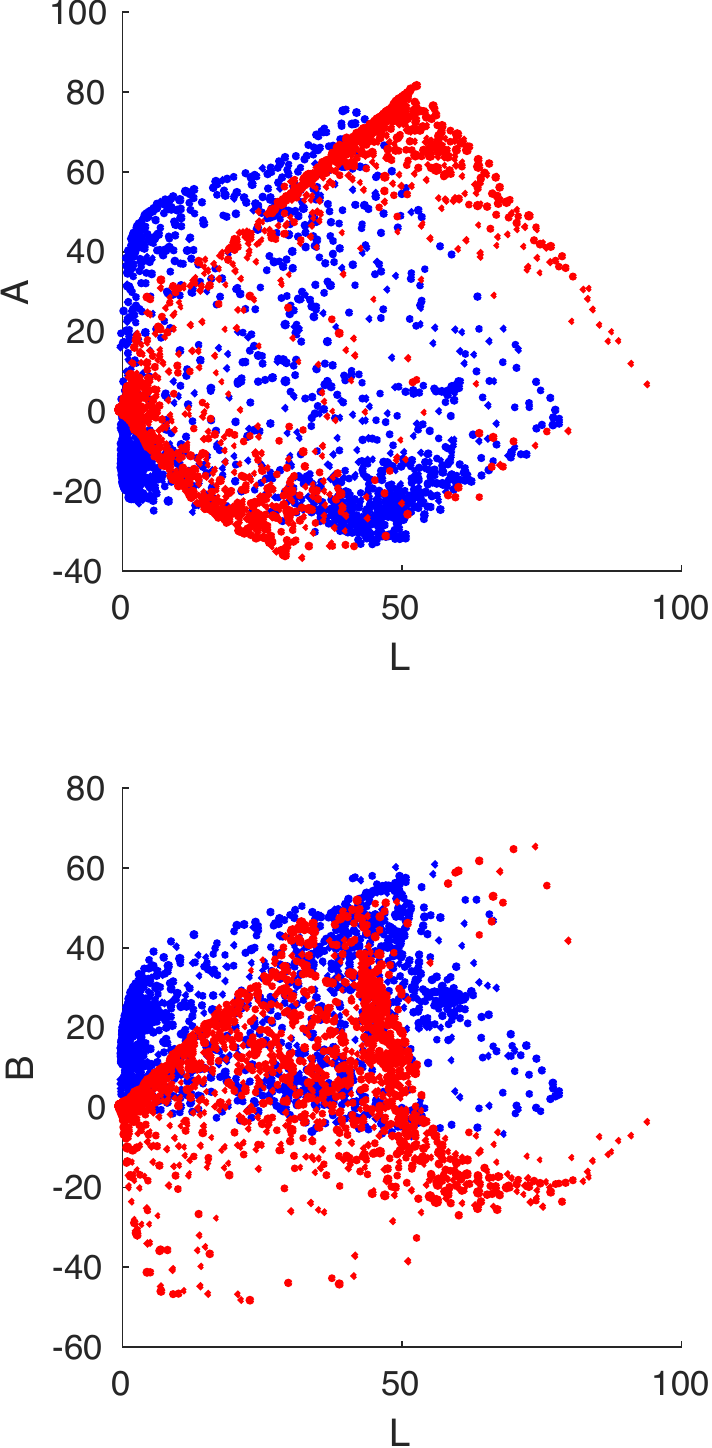}}
\subfloat[COT result in LAB]{\includegraphics[width = 0.35\textwidth]{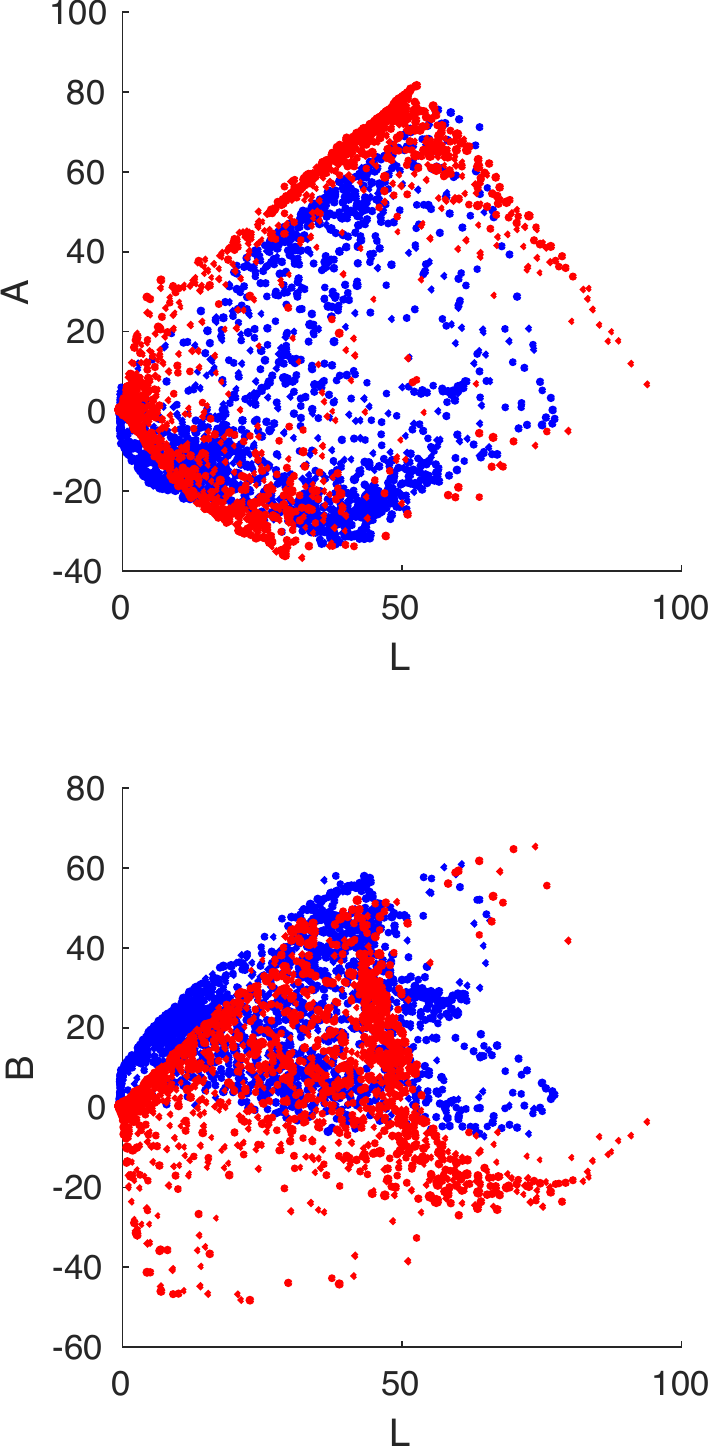}}
\subfloat[3DOT result in LAB]{\includegraphics[width = 0.35\textwidth]{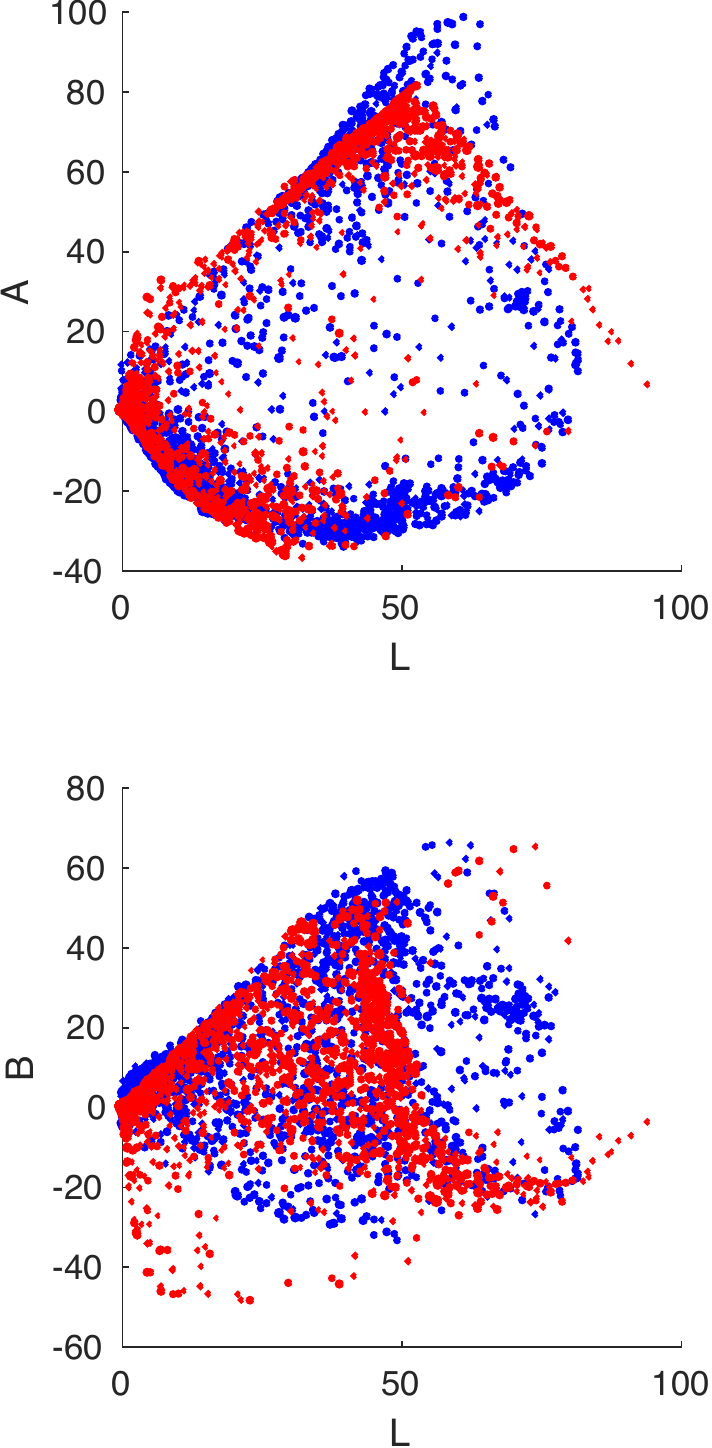}}\\
\caption{Left Column: image obtained performing one dimensional optimal transport for the Luminosity (L) coordinate ignoring the A and B coordinates. Second column: image obtained by performing optimal transport on luminosity conditioned on color. Third column: plain three dimensional optimal transport in L*a*b space.}
\label{fig:red_flowers}
\end{figure}

Figure  \ref{fig:red_flowers} shows the result obtained changing lightness in three different ways. First (left column) we use one-dimensional optimal transport (with quadratic cost) to map the $L$ coordinate, ignoring the values of $A$ and $B$. The L*a*b diagram shows that this results  in a nearly uniform shift of $L$ towards smaller values. The third column shows the effect of mapping the starting image to the target image through 3d optimal transport in the full L*a*b space. In this case the point clouds overlap to a much better degree, yet we observe that the color of the lotus has been changed too much towards the color on the poinsettia of the target image. The second column is obtained performing optimal transport of $L$ conditioned on $A$ and $B$. Contrasting to the other two results, here the lotus has kept its original color, and the lightness has changed to a different degree for the lotus than for the background leaves.

This is a general advantage of conditional optimal transport: unlike its unconditional cousin, it does not need to preserve total mass (in this case, transferring fully one color palette to the other), but only the mass for each value of $z$. This point to an additional application of conditional optimal transport: its capacity to address possible unbalances between source and target by parameterizing the transfer map by means of convenient labels $z$. In work in progress, we expand on this notion, finding those latent covariates $z$ that help resolve unbalances optimally.

\newpage

\bibliographystyle{plain}      
\bibliography{GB}

\end{document}